\documentclass[3p]{elsarticle}

\usepackage{verbatim}
\usepackage[pdftex]{color}
\usepackage{amssymb}
\usepackage[colorlinks,pdfpagelabels,pdfstartview = FitH,bookmarksopen = true,bookmarksnumbered = true,linkcolor = blue,plainpages = false,hypertexnames = false,citecolor = red] {hyperref}
\usepackage[active]{srcltx}
\usepackage[final]{pdfpages} 

\newtheorem{rmk}{Remark}

\usepackage{amsfonts,euscript,eufrak,latexsym,amssymb}
\usepackage{amsmath, amssymb, graphicx}
\usepackage[english]{babel}
\usepackage[latin1]{inputenc}
\usepackage[T1]{fontenc}

\usepackage{mathtools}

\newtoks\by
\newtoks\paper
\newtoks\book
\newtoks\jour
\newtoks\yr
\newtoks\pages
\newtoks\vol
\newtoks\publ

\def\name[#1, #2]{#1 #2}
\def\ota{{\hbox{\bf ???}}}
\def\cLear{\by=\ota\paper=\ota\book=\ota\jour=\ota\yr=\ota
\pages=\ota\vol=\ota\publ=\ota}
\def\endpaper{\the\by, \textit{\the\paper},
{\the\jour} \textbf{\the\vol} (\the\yr), \the\pages.\cLear}
\def\endbook{\the\by, \textit{\the\book},
\the\publ, \the\yr.\cLear}
\def\endpap{\the\by, \textit{\the\paper}, \the\jour.\cLear}
\def\endproc{\the\by, \textit{\the\paper}, \the\book, \the\publ,
\the\yr, \the\pages.\cLear}

\begin{document}

\title{Nonlinear Reynolds equation for hydrodynamic lubrication}
\author[aalto]{Tom Gustafsson}
\ead{tom.gustafsson@aalto.fi}
\author[tamu]{K.R. Rajagopal\corref{cor1}}
\ead{krajagopal@tamu.edu}
\author[aalto]{Rolf Stenberg}
\ead{rolf.stenberg@aalto.fi}
\author[ist]{Juha Videman}
\ead{videman@math.ist.utl.pt}

\address[aalto]{Aalto University,
Department of Mathematics and Systems Analysis, \\ P.O. Box 11100
FI-00076 Aalto,
Finland}

\address[tamu]{Department of Mechanical Engineering, Texas A\&M University, \\ 3123 College Station, TX 77843-3123, USA}

\address[ist]{CAMGSD and
Mathematics Department,
Instituto Superior T\'ecnico, \\
Universidade de Lisboa, 
Av. Rovisco Pais 1, 1049-001 Lisboa, Portugal}

\cortext[cor1]{Corresponding author. Tel. +1-521-820-0782}

\begin{abstract}
We derive a novel and rigorous  correction to the classical Reynolds lubrication approximation  for fluids with viscosity depending upon the pressure. Our analysis shows that the pressure dependence of viscosity leads to additional nonlinear terms related to the shear-rate and arising from a non negligible cross-film pressure. We present a numerical comparison between the classical Reynolds equation and our modified equation and conclude that the  modified equation  leads to the prediction of higher pressures and viscosities in the flow domain.
\end{abstract}

\maketitle


\section{Introduction}

The Reynolds equation \cite{Rey1886}, which is an approximation of the classical Navier--Stokes equations, describes reasonably well the flow of a large class of fluids whose viscosities can be assumed to be independent of the pressure in many lubrication problems. However, there is a clear and incontestable evidence, that is well-documented, that attests to the fact that the viscosity varies with pressure, especially so in several problems concerning thin film lubrication, cf.~\cite{Bridgman31,Bair2007,Szeri2011}. Experiments have shown that in high pressure regimes  pressure variations can  significantly change the viscosity of certain lubricants and in areas such as elastohydrodynamic lubrication the classical isoviscous lubrication theory is noticeably incapable of explaining the existence of continuous lubricant films, for example, in rolling-contact bearings, cf.~\cite{Szeri2011}. 

The traditional correction to the Reynolds approximation in the piezoviscous regime is based on a rather heuristic assumption that it suffices to replace the constant viscosity in the Reynolds equation by a suitable viscosity-pressure relationship, cf.~\cite{DH66}.  This approach becomes  questionable, however, in high pressure regimes because of the possible change of type  (loss of ellipticity) of the equations and the potential existence of cross-film pressure gradient, see the discussion in \cite {RS03}. In fact, for a Reynolds type approximation to be valid in this regime it ought be  derived from the full balance of linear momentum equations governing the flow of incompressible fluids with pressure-dependent viscosities. Now,  the mathematical theory for the equations that govern the flows of fluids with a pressure-dependent viscosity has been advancing in leaps  and bounds over the last few decades, see  \cite{Ren86,Gazz97,HMR01,MNR02,FMR05,MR07,BMR09,SV10,HLS12}, but the only rigorous attempt to derive a modified Reynolds equation for elastohydrodynamic lubrication based upon these equations seems to have been carried out by Rajagopal and Szeri \cite{RS03}, see also \cite{BCGV13} for further applications of their model.

In this paper we propose a new modified Reynolds equation for hydrodynamic lubrication in high pressure regimes. Our approach is built upon an asymptotic expansion of the non-dimensional velocity and pressure fields in terms of a small dimensionless parameter $\epsilon$, related to the film thickness, and on a systematic analysis of the simplified set of equations obtained at different orders of $\epsilon$, see \cite{NV07} for a similar approach in the isoviscous case.  Assuming  that the dimensionless pressure-viscosity coefficient is of order $\epsilon$, we show that the pressure distribution is governed by a Reynolds equation modified by a term depending, in particular, on the square of the shear rate. In \cite{RS03}, the authors assumed, for simplicity, that the cross-film pressure vanishes and derived a modified Reynolds equation, similar to ours, but depending on the elongation rate. Our computations show that  it is exactly the (lower-order) cross-film pressure which is responsible for the new modified term in the Reynolds equation. 
We also prove that for smaller values of the pressure-viscosity coefficient the traditional modification of the Reynolds equation is a very accurate approximation of the pressure field. 

It has been argued that the behavior of a piezoviscous fluid cannot be adequately described by a Reynolds type equation  if the principal shear stress $\tau$  is not  less than the reciprocal of the pressure-viscosity coefficient $\alpha$, cf.~\cite{SGRW00}. This is not surprising since problems of nonexistence and nonuniqueness are expected for the full balance of linear momentum  equations if $\tau\geq (\alpha)^{-1}$, cf.~\cite{Ren86}. Similar conclusions can also be drawn from our modified Reynolds equations which ceases to be elliptic if the shear stress times the pressure-viscosity coefficient is larger than, or of the order of, one.

% After laying down the initial assumptions leading to our lubrication approximation we do not  make any  further simplifying assumptions on the while deriving  a modified Reynolds equation, nor do we simplify our modified equation for numerical computations. 

In the piezoviscous regime,  the pressure-dependence of viscosity can be modeled through the Barus relation (Barus \cite{Barus93})
\begin{equation}
\mu=\mu_0\, \mathrm{e}^{\alpha  p}\, ,
\label{barus}
\end{equation}
 where $\mu_0$  denotes the dynamic viscosity measured at the ambient pressure ($p=0$) and $\alpha$ is a positive parameter related to the rate of change of viscosity with respect to pressure, often referred to as the pressure-viscosity coefficient.  The  Barus relationship,  together with the Roelands formula (Roelands \cite{Roe66})
 \begin{equation}
 \mu=\mu_0 \, \left( \frac{\mu_0}{\mu_R}\right)^{1-(1+p/p_R)^{Z}} \, ,
 \label{roelands}
 \end{equation}
 where $\mu_R=6.31\cdot 10^{-5}$ Pa s, $p_R=1.98\cdot 10^8$ Pa and $Z$ is a dimensionless parameter usually adjusted at ambient pressure to the Barus relation through the equation
 \begin{equation}
 Z=\frac{\alpha \, p_R}{\ln \mu_0-\ln \mu_R} \, ,
 \label{zform}
 \end{equation}
are the most widely used models in elastohydrodynamic lubrication. For applications of these and other experimentally validated formulas for the variation of viscosity with pressure, temperature and density, see, e.g., \cite{Bridgman31,PDR99,BK03,BLW06,Bair2007,Szeri2011}.

In general, the pressure-viscosity coefficient $\alpha$ depends on the lubricant, and on the pressure, temperature and shear rate in the contact area, cf.~\cite{Bair2007}. For different lubricant oils, its value has  been shown to vary between $1 \cdot 10^{-8}\,  ({\rm Pa})^{-1}$ and $4\cdot 10^{-8}\, ({\rm Pa})^{-1}$, cf.~\cite{HSJ04,vanLee09,PSR12}. If the constant $Z$ in the Roelands formula \eqref{roelands} is computed from equation \eqref{zform}, then both  the Barus and the Roelands relation give $\frac{1}{\mu}\frac{{\rm d} \mu}{{\rm d} p} =\alpha$. We stress that although we  rely on the assumption that the viscosity $\mu$ varies with pressure $p$ according to the Barus or Roelands  formula, our computations hold for more general pressure-viscosity relationships for which $\frac{1}{\mu}\frac{{\rm d} \mu}{{\rm d} p}$ is of the order of $\epsilon$ or smaller.

Before concluding this introductory section, it is worthwhile recognizing a marked departure of the constitutive relation  of a piezoviscous fluid from that of the classical incompressible Navier--Stokes fluid with constant viscosity. While the latter is described by an explicit constitutive function for the stress in terms of the symmetric part of the velocity gradient, the former is an implicit relationship for the Cauchy stress and the symmetric part of the velocity gradient, leading to a totally different structure to the equations governing the flows of such fluids which in turn raises interesting issues with regard to both the mathematical and numerical analysis of the governing equations.

A one-dimensional rate-type implicit model to describe the non-Newtonian response of fluids was introduced by Burgers \cite{Bur39}. His model includes as a special case the pioneering model to describe the viscoelastic response of fluids that was advanced by Maxwell \cite{Max66}. Maxwell's model is however not an implicit model as the symmetric part of the velocity gradient can be expressed explicitly in terms of the stress and the time rate of the stress. Oldroyd \cite{Old50} developed a systematic procedure to generate properly invariant three-dimensional rate type implicit constitutive relations. Such fluid models can be used to describe the flow of viscoelastic fluids and it is our aim to generalize the type of approximation that is being carried out here to include rate type implicit constitutive relations.

\section{Lubrication approximation}

Consider the following equations governing the isothermal flow of an  incompressible, homogeneous, viscous fluid
\begin{align}
\displaystyle \rho\, \left( \frac{\partial \mathbf{v}}{\partial t} + (\mathbf{v}\cdot\nabla)\, \mathbf{v} \right) -2\, \nabla\cdot \big(\, \mu(p)  \mathbf{D}(\mathbf{v})\,\big) +\nabla p & = \rho \mathbf{b \, ,}
\label{nseqs1} \\  \nonumber \\ 
\nabla\cdot \mathbf{v} & =0 \, ,
\label{nseqs2}
\end{align} 
where $\rho>0$ is the constant density of the fluid, $\mathbf{v}= (u,v,w)$ is the velocity field and $p$ is a scalar variable, often referred to as the mechanical pressure, associated with the incompressibility constraint \eqref{nseqs2}. Moreover, we have assumed above that the stress response is linear in the symmetric part of the velocity gradient, $\mathbf{D}(\mathbf{v})=\frac 12\left(\nabla \mathbf{v}+(\nabla\mathbf{v})^T\right)$, with the viscosity $\mu$ depending (in isothermal conditions) only on the pressure $p$, so that the stress tensor $\mathbf{T}$ reads as
\begin{equation}
\mathbf{T}=-p\mathbf{I}+ 2\mu(p)\, \mathbf{D}(\mathbf{v})\, .
\label{implicitlaw}
\end{equation}
Ignoring the body force  $\mathbf{b}$ and limiting oneself to steady motions, equations \eqref{nseqs1}--\eqref{nseqs2}  become
\begin{align}
\displaystyle \rho\, (\mathbf{v}\cdot\nabla)\, \mathbf{v} -  \mu(p)\Delta \mathbf{v} -2 \mathbf{D}(\mathbf{v})\, \nabla \mu(p)  +\nabla p & = 0\,,  \label{steadyns1}  \\  \nonumber \\ 
\nabla\cdot \mathbf{v} & =0 \, .
\label{steadyns2}
\end{align}
Note that \eqref{implicitlaw} defines an implicit relationship $\mathbf{f}(\mathbf{T},\mathbf{D})=\mathbf{0}$ between the stress tensor $\mathbf{T}$ and the symmetric part of the velocity gradient  $\mathbf{D}$ since  in an incompressible fluid ${\rm tr}\, \mathbf{D}=0$ so that  $p=-\, \frac{1}{3} {\rm tr}\, \mathbf{T}$, i.e.~the mechanical pressure is  just the mean normal stress. 

Let us assume, for expediency, that the viscosity depends on the pressure through the Barus relation \eqref{barus}.
Introducing the non-dimensionalized (starred) quantities 
\[
\mathbf{x}^\ast= L^{-1} \mathbf{x}\, , \quad \mathbf{v}^\ast = U^{-1} \mathbf{v} \, , \quad  \mu^\ast=\mu_0^{-1}\mu\, ,  \quad p^\ast=P^{-1}p \, ,  \quad \alpha^\ast=P\, \alpha\, ,
\]
where $L$ and $U$ represent typical length and velocity scales and the characteristic pressure is taken to be
$
P=\mu_0 U L^{-1} \, ,
$
 we can rewrite equations \eqref{steadyns1}--\eqref{steadyns2}  as
\begin{align}
\mu^\ast\, (p^\ast) \Delta^\ast \mathbf{v}^\ast +2\, \alpha^\ast\mu^\ast(p^\ast) \, \mathbf{D}^\ast(\mathbf{v}^\ast)\, \nabla^\ast p^\ast -\nabla^\ast p^\ast & =\mathrm{Re}\, (\mathbf{v}^\ast\cdot\nabla^\ast)\mathbf{v}^\ast  \, , \label{ndeqs1}
\\ \nonumber \\ 
\nabla^\ast\cdot \mathbf{v}^\ast & =0\, ,
\label{ndeqs2}
\end{align}
where $\mathrm{Re}=\rho\, UL\mu_0^{-1}$ denotes the usual Reynolds number. Note that, had we defined $\alpha$ more generally through
\begin{equation}
\alpha=\frac{1}{\mu}\frac{{\rm d} \mu}{{\rm d} p}\, ,
\label{pvform}
\end{equation}
the form of the  nondimensional system \eqref{ndeqs1}--\eqref{ndeqs2} would still remain the same.

Let us restrict our attention to two-dimensional plane flows
\[
\mathbf{v^\ast}(\mathbf{x^\ast})=\big(u^\ast(x^\ast,y^\ast), v^\ast(x^\ast,y^\ast)\big)\, , \qquad p^\ast=p^\ast(x^\ast,y^\ast)\, ,
\]
so that equations \eqref{ndeqs1}--\eqref{ndeqs2} can be recast  as
\begin{align}
\displaystyle \mu^\ast  \left(\frac{\partial^2 u^\ast }{\partial {x^\ast}^2}  + \frac{\partial^2 u^\ast }{\partial {y^\ast}^2}\right)  & + \alpha^\ast\mu^\ast\left(2\, \frac{\partial u^\ast}{\partial x^\ast}  \frac{\partial p^\ast}{\partial x^\ast}  +
\left(\frac{\partial u^\ast}{\partial y^\ast} +\frac{\partial v^\ast}{\partial x^\ast} \right) \frac{\partial p^\ast}{\partial y^\ast} \right) -
\frac{\partial p^\ast}{\partial x^\ast}  \nonumber \\  \nonumber\\ 
\displaystyle    \qquad \qquad \quad & =\mathrm{Re}\, \left( u^\ast\frac{\partial u^\ast}{\partial x^\ast} +v^\ast \frac{\partial u^\ast}{\partial y^\ast}  \right)\, , \label{twodeqs1} \\  \nonumber \\ 
\displaystyle \mu^\ast  \left(\frac{\partial^2 v^\ast }{\partial {x^\ast}^2}  + \frac{\partial^2 v^\ast }{\partial {y^\ast}^2}\right) & + \alpha^\ast\mu^\ast\left( 
\left(\frac{\partial u^\ast}{\partial y^\ast} +\frac{\partial v^\ast}{\partial x^\ast} \right) \frac{\partial p^\ast}{\partial x^\ast}  + 2\, \frac{\partial v^\ast}{\partial y^\ast}  \frac{\partial p^\ast}{\partial y^\ast}  \right) -
\frac{\partial p^\ast}{\partial y^\ast}  \nonumber  \\  \nonumber \\ 
\displaystyle  \qquad \qquad \quad  & =\mathrm{Re}\, \left( u^\ast\frac{\partial v^\ast}{\partial x^\ast} +v^\ast \frac{\partial v^\ast}{\partial y^\ast}  \right) \, ,
\label{twodeqs2}\\  \nonumber \\ 
\displaystyle \frac{\partial u^\ast}{\partial x^\ast} +\frac{\partial v^\ast}{\partial y^\ast} & = 0 \, .
\label{twodeqs3}
\end{align}
We make the following  assumptions

\begin{itemize}

\item the flow takes place between two almost parallel surfaces situated at $y^\ast=0$ and $y^\ast=H^\ast(x^\ast)$; 

\item curvature of the surfaces can be neglected;

\item the lubrication film is thin,  that is, $H^\ast(x^\ast)=\epsilon h^\ast(x^\ast)$, where $\epsilon \ll 1 $ denotes a small non-dimensional parameter;

\item the characteristic lengths in $x$- and $y$-directions scale as $L_y/L_x=\mathcal{O}(\epsilon)$; 

\item the flow is slow enough or the viscosity high enough so that $\mathrm{Re}=\mathcal{O}(\epsilon)$;

\item the scaled viscosity is of the order of one, that is $\mu=\mathcal{O}(1)$. 

\begin{rmk}
Choosing $\mathrm{Re}=\mathcal{O}(\epsilon)$ we end up neglecting the inertial effects from the outset, as usual in deriving the Reynolds approximation. At the same time, we simplify our presentation. In fact,  given that the main pressure approximation is not affected by the inertial effects for Reynolds numbers up to the order $\mathcal{O}(\epsilon^{-1})$, cf.~\cite{NV07} for  computations in the constant-viscosity case,    this simplification seems entirely warranted.  

As for the assumption $\mu=\mathcal{O}(1)$, it is essential for the dimension reduction and, therefore, for the existence of a Reynolds type equation. It may not hold at the highest pressure regimes but, then again, for pressures higher than 0.5 GPa the entire viscosity-pressure relationships given by the Barus or Roelands formula becomes questionable, cf.~\cite{Bair2007, vanLee09}. 
\end{rmk}

\end{itemize}
Redefining the $y$-variable as $y=\epsilon^{-1}y^\ast$ and dropping the stars in \eqref{twodeqs1}--\eqref{twodeqs3}, yields the system
\begin{align}
    \displaystyle \mu  \left(\frac{\partial^2 u }{\partial x^2} +\epsilon^{-2} \frac{\partial^2 u }{\partial y^2}\right)  &+  \alpha \mu\left(2\, \frac{\partial u}{\partial x}  \frac{\partial p}{\partial x}  + \epsilon^{-2} \frac{\partial u}{\partial y}\, \frac{\partial p}{\partial y}  +\epsilon^{-1} \, \frac{\partial v}{\partial x} \, \frac{\partial p}{\partial y} \right) -
\frac{\partial p}{\partial x} \nonumber \\\nonumber \\ 
\displaystyle &=\mathrm{Re}\, \left( u\frac{\partial u}{\partial x} +\epsilon^{-1}\, v \frac{\partial u}{\partial y}  \right)\, ,  \label{epseqs1} \\\nonumber  \\ \nonumber
\displaystyle \mu  \left(\frac{\partial^2 v }{\partial x^2} +\epsilon^{-2} \frac{\partial^2 v }{\partial y^2}\right)  &+ \alpha\mu\left( \epsilon^{-1}
\frac{\partial u}{\partial y} \, \frac{\partial p}{\partial x} + \frac{\partial v}{\partial x} \, \frac{\partial p}{\partial x}  + 2 \epsilon^{-2} \frac{\partial v}{\partial y}  \frac{\partial p}{\partial y}  \right) -
\epsilon^{-1} \frac{\partial p}{\partial y}  \label{epseqs2}  \\ \nonumber \\ 
\displaystyle   &=\mathrm{Re}\, \left(  u\frac{\partial v}{\partial x} +\epsilon^{-1}\,v\frac{\partial v}{\partial y}  \right) \, ,
\\ \nonumber \\ 
\displaystyle \frac{\partial u}{\partial x} +\epsilon^{-1} \frac{\partial v}{\partial y}  &= 0 \, ,  \label{epseqs3} \\  \nonumber  \\ 
\mathbf{v}=\mathbf{U}^0\, , \quad \mathbf{U}^0 \cdot \mathbf{n}&=0,\qquad \mathrm{at} \ \ y=0\,,   \label{epseqs4} \\  \nonumber  \\ 
\mathbf{v}=\mathbf{U}^h\, , \quad \mathbf{U}^h\cdot \mathbf{n}&=0,\qquad \mathrm{at} \ \ y=h(x)\, ,
\label{epseqs5}
\end{align}
 where $\mathbf{U}^0$ and $\mathbf{U}^h$ denote the  velocities of the parallel plates. We will write 
 \[
 \mathrm{Re}=\epsilon \, \texttt{Re}\]
  and assume that
\begin{align}
\mathbf{v}(\epsilon,\mathbf{x})=\epsilon^0\mathbf{v}^1(\mathbf{x})+\epsilon^1\mathbf{v}^2(\mathbf{x})+\epsilon^2\mathbf{v}^3(\mathbf{x})+\ldots\, , \label{ansatz1}
\\  \nonumber \\ 
p(\epsilon,\mathbf{x})=\epsilon^{-2}p^0(x)+\epsilon^{-1} p^1(\mathbf{x})+\epsilon^0p^2(\mathbf{x})+\ldots\, ,
\label{ansatz2}
\end{align}
where  the functions $\mathbf{v}^j$ and $p^j$ and the parameter $\texttt{Re}$ are of  $\mathcal{O}(1)$. 

There is still one dimensionless parameter, $\alpha$, whose order has not been discussed. In the thinnest-film (elastohydrodynamic)  lubrication problems, $\epsilon$ is of the order $10^{-6}$ or smaller, and the size of $\alpha$ should be taken to be of the order of $\epsilon$ but in thicker film problems the order of the pressure-viscosity coefficient is $\epsilon^2$ or smaller. We will consider these two cases separately. 

\subsection{Case $\alpha=\mathcal{O}(\epsilon)$}

Substituting the asymptotic expressions \eqref{ansatz1}--\eqref{ansatz2}  in equations \eqref{epseqs1}--\eqref{epseqs5}, writing 
\[
\alpha=\epsilon\, \mathcal{A}\,,
\]
where $ \mathcal{A}=\mathcal{O}(1)$, and keeping only the terms of the highest order, we obtain, at order $\epsilon^{-2}$ in \eqref{epseqs1} and \eqref{epseqs2},  at order $\epsilon^{-1}$ in \eqref{epseqs3} and at order  $\epsilon^0$ in \eqref{epseqs4} and  \eqref{epseqs5}
\begin{align}
    \displaystyle \mu  \frac{\partial^2 u^1 }{\partial y^2} + \mathcal{A}\mu  \frac{\partial u^1}{\partial y}  \frac{\partial p^1}{\partial y} - \frac{\partial p^0}{\partial x} &= 0 \, , \label{fstorder1}\\ \nonumber \\
\displaystyle \mu  \frac{\partial^2 v^1 }{\partial y^2} + \mathcal{A} \mu  \left(\frac{\partial u^1}{\partial y}  \frac{\partial p^0}{\partial x}+ 2\, \frac{\partial v^1}{\partial y}  \frac{\partial p^1}{\partial y} \right) - \frac{\partial p^1}{\partial y} &= 0 \, , \label{fstorder2}
\\ \nonumber \\
\frac{\partial v^1}{\partial y}  &= 0\,, \label{fstorder3} \\ \nonumber \\
u^1=U^0\,, \quad v^1&=0\, , \qquad   \mathrm{at} \ \ y=0 \, , \label{fstorder4} \\\nonumber  \\ 
u^1=U^h\, , \quad v^1&=0\, , \qquad   \mathrm{at} \ \ y=h(x) \, . \label{fstorder5} 
\end{align}
Equations \eqref{fstorder3}, \eqref{fstorder4} and \eqref{fstorder5}, show that
\[
v^1(\mathbf{x})\equiv  0\, .
\]
It then follows from  \eqref{fstorder1} and \eqref{fstorder2} that
\begin{equation}
\mu \frac{\partial^2 u^1 }{\partial y^2}  +  \mathcal{A}^2\mu^2 \, \frac{\partial p^0}{\partial x}  \left( \frac{\partial u^1}{\partial y} \right)^2 \, =  \frac{\partial p^0}{\partial x}  \, . \label{reyn1}
\end{equation}

\bigskip

The next order equations read as
\begin{align}
    \displaystyle \mu  \frac{\partial^2 u^2 }{\partial y^2} +\mathcal{A}  \mu \frac{\partial p^1}{\partial y}   \frac{\partial u^2}{\partial y} &= \frac{\partial p^1}{\partial x} - 2 \mathcal{A} \mu  \frac{\partial u^1}{\partial x}  \frac{\partial p^0}{\partial x}  -  \mathcal{A} \mu  \frac{\partial u^1}{\partial y}  \frac{\partial p^2}{\partial y}    \, , \label{secorder1}\\ \nonumber \\
    \displaystyle \mu  \frac{\partial^2 v^2 }{\partial y^2}  +2\, \mathcal{A} \mu \frac{\partial p^1}{\partial y}   \frac{\partial v^2}{\partial y} - \frac{\partial p^2}{\partial y} &= - \mathcal{A} \mu  \frac{\partial u^1}{\partial y}  \frac{\partial p^1}{\partial x}- \mathcal{A} \mu \frac{\partial p^0}{\partial x}   \frac{\partial u^2}{\partial y}  \, , \label{secorder2}
\\ \nonumber \\
\frac{\partial v^2}{\partial y} &= -\frac{\partial u^1}{\partial x} \,,  \label{secorder3} \\ \nonumber \\
u^2=v^2&=0\, ,  \qquad \mathrm{at} \ \ y=0 \, , \label{secorder4} \\\nonumber  \\ 
v^2=h^\prime U^h\,, \quad u^2 &= 0 \, ,    \qquad \mathrm{at} \ \ y=h(x) \, . \label{secorder5} 
\end{align}
The one-dimensional Oseen problem \eqref{secorder2}--\eqref{secorder5} for $(v^2,p^2)$ is solvable if and only if the compatibility condition
\[
-\, \int_0^{h} \frac{\partial u^1}{\partial x} \, \mathrm{d} y = h^\prime\,  U^h \, 
\]
is satisfied. This condition can be written as
\begin{equation}
 \frac{\mathrm{d} }{\mathrm{d} x}  \,  \int_0^{h} u^1 \, \mathrm{d} y =0 \, .
 \label{comp}
 \end{equation}
 Multiplying \eqref{reyn1} by $y(y-h)$ and integrating across the film,  yields 
\begin{equation}
2 \mu\,   \int_0^{h} u^1 \, \mathrm{d} y +\mu h(U^0-U^h) +\mathcal{A}^2\mu^2 \, \frac{\partial p^0}{\partial x} \, \int_0^h y(y-h)\, \left( \frac{\partial u^1}{\partial y} \right)^2 \, \mathrm{d} y = -\, \frac{h^3}{6 }\frac{\partial p^0}{\partial x} \, .
 \label{calc}
 \end{equation}
Using \eqref{comp} in \eqref{calc}  leads to  the modified Reynolds equation
\[
\frac{\mathrm{d}}{\mathrm{d} x} \left[  \left( \frac{h^3}{12\mu }  - \mathcal{A}^2\mu  \, \int_0^h y(h-y)\, \left( \frac{\partial u^1}{\partial y} \right)^2 \, \mathrm{d} y \right) 
 \frac{\mathrm{d} p^0}{\mathrm{d} x}    \right] = \frac{h^\prime}{2}(U^h-U^0)\, .
\]

\subsection{Case $\alpha=\mathcal{O}(\epsilon^2)$}

At the first order, we simply have
\begin{align}
    \displaystyle \mu  \frac{\partial^2 u^1 }{\partial y^2}  &= \frac{\partial p^0}{\partial x}  \, , \label{forder1} \\ \nonumber \\
    \displaystyle \mu  \frac{\partial^2 v^1 }{\partial y^2} - \frac{\partial p^1}{\partial y} &= 0 \, , \label{forder2}
\\ \nonumber \\
\frac{\partial v^1}{\partial y} &= 0\,, \label{forder3} \\ \nonumber  \\
u^1=U^0\,, \quad v^1&=0\, , \qquad   \mathrm{at} \ \ y=0 \, , \label{forder4}  \\ \nonumber \\
u^1=U^h\, , \quad v^1&=0\, , \qquad   \mathrm{at} \ \ y=h(x) \, . \label{forder5} 
\end{align}
From \eqref{forder3}--\eqref{forder5} it again follows that $v^1(\mathbf{x})\equiv 0$.  Equation \eqref{forder2} then yields $\partial p^1/\partial y=0$, i.e.~the  cross-film pressure also vanishes  at this order. Therefore 
\begin{equation}
u^1(x,y)= U^0+\frac{U^0-U^h}{h} \, y + \frac{y(y-h)}{2\mu}\, \frac{\partial p^0}{\partial x} \, .
\label{uexpr}
\end{equation}
At the next order, we find that
\begin{align}
    \displaystyle \mu  \frac{\partial^2 u^2 }{\partial y^2} &= \frac{\partial p^1}{\partial x}  \, , \label{sorder1}\\ \nonumber \\
    \displaystyle \mu  \frac{\partial^2 v^2 }{\partial y^2}  - \frac{\partial p^2}{\partial y} &= - \mathcal{A} \mu  \frac{\partial u^1}{\partial y}  \frac{\partial p^0}{\partial x}  \, , \label{sorder2}
\\ \nonumber \\
\frac{\partial v^2}{\partial y}  &= -\frac{\partial u^1}{\partial x} \,,  \label{sorder3} \\ \nonumber \\
u^2=v^2&=0\, , \qquad   \mathrm{at} \ \ y=0 \, , \label{sorder4} \\\nonumber  \\ 
v^2=h^\prime U^h\, , \quad u^2&=0\,, \qquad   \mathrm{at} \ \ y=h(x) \, . \label{sorder5} 
\end{align}
The one-dimensional Stokes problem \eqref{sorder2}--\eqref{sorder5} for $(v^2,p^2)$ is solvable if and only if 
\begin{equation}
-\, \int_0^{h} \frac{\partial u^1}{\partial x} \, \mathrm{d} y = h^\prime\,  U^h \, .
\label{compat}
\end{equation}
Substituting \eqref{uexpr} into \eqref{compat}, yields
\[
-\, \frac{\partial}{\partial x} \, \left( hU^0-(U^0-U^h)\frac{h^2}{2} +\frac{1}{2\mu}\, \int_0^h y(y-h)  \, \mathrm{d} y  \,  \frac{\partial p^0}{\partial x}\right) + h^\prime \, U^h =h^\prime \, U^h\, ,
\]
which can be written as the  classical Reynolds equation 
\begin{equation}
\frac{1}{12}\, \frac{\mathrm{d}}{\mathrm{d} x} \, \left( \frac{h^3}{\mu}  \frac{\mathrm{d} p^0}{\mathrm{d} x}\right) = \frac{h^\prime}{2}\left( U^h+U^0\right)\, .
\label{reynolds}
\end{equation}

\bigskip

At the next order, we obtain the following system for $(\mathbf{v}^3,p^3)$ 
\begin{align}
    \displaystyle \mu  \frac{\partial^2 u^3 }{\partial y^2} &= \frac{\partial p^2}{\partial x} - 2 \mathcal{A} \mu  \frac{\partial u^1}{\partial x}  \frac{\partial p^0}{\partial x}  -  \mathcal{A} \mu  \frac{\partial u^1}{\partial y}  \frac{\partial p^2}{\partial y}    \, , \label{torder1}\\ \nonumber \\
    \displaystyle \mu  \frac{\partial^2 v^3 }{\partial y^2}  - \frac{\partial p^3}{\partial y} &= - \mathcal{A} \mu  \frac{\partial u^1}{\partial y}  \frac{\partial p^1}{\partial x} \, , \label{torder2}
\\ \nonumber \\
\frac{\partial v^3}{\partial y}  &= -\frac{\partial u^2}{\partial x} \,,  \label{torder3} \\ \nonumber \\
u^3=v^3&=0\, , \qquad   \mathrm{at} \ \ y=0 \, , \label{torder4} \\\nonumber  \\ 
u^3=-\frac{1}{2}{h^\prime}^2 U^h\, , \quad v^3&=0\, , \qquad   \mathrm{at} \ \ y=h(x) \, . \label{torder5} 
\end{align}
The one-dimensional Stokes problem \eqref{torder2}--\eqref{torder5} is solvable if and only if 
\[
-\, \int_0^{h} \frac{\partial u^2}{\partial x} \, \mathrm{d} y = 0\, .
\]
On the other hand, from \eqref{sorder1},  \eqref{sorder4} and  \eqref{sorder5}  one concludes that
\[
u^2(x,y)=  \frac{y(y-h(x))}{2\mu}\, \frac{\partial p^1}{\partial x} \,.
\label{u2expr}
\]
This yields
\[
\frac{1}{12}\, \frac{\mathrm{d}}{\mathrm{d} x}\left( \frac{h^3}{ \mu} \, \frac{\mathrm{d} p^1}{\mathrm{d} x}\right) =  0\, .
\]
We thus conclude that the Reynolds equation \eqref{reynolds} provides a very accurate approximation for the pressure when  $\alpha=\mathcal{O}(\epsilon^2)$ or smaller.

\section{Numerical results}

Contrary to \cite{RS03}, we do not  simplify our modified Reynolds equation for numerical computations. Rather, we couple it with a nonlinear ODE governing the behavior of the along channel velocity 
and solve simultaneously for the pressure and the main (along-channel) velocity component. 
The modified Reynolds equation 
\begin{equation}
\frac{\mathrm{d}}{\mathrm{d} x} \left[  \left( \frac{h^3}{12\mu }  - \mathcal{A}^2\mu  \, \int_0^h y(h-y)\, \left( \frac{\partial u^1}{\partial y} \right)^2 \, \mathrm{d} y \right) 
 \frac{\mathrm{d} p^0}{\mathrm{d} x}    \right] = \frac{h^\prime}{2}(U^h-U^0)\,
\label{modified}
\end{equation}
is a nonlinear second-order ODE for $p^0$. Complemented with suitable boundary conditions, it becomes  solvable when considered with  equation
\begin{align}
    \mu \frac{\partial^2 u^1 }{\partial y^2}  +  \mathcal{A}^2\mu^2 \, \frac{\partial p^0}{\partial x}  \left( \frac{\partial u^1}{\partial y} \right)^2 \, &=  \frac{\partial p^0}{\partial x}  \, , \label{usol1} \\ \nonumber \\ 
    u^1&=U^0\, , \qquad   \mathrm{at} \ \ y=0 \, , \label{usol2} \\\nonumber  \\ 
    u^1&=U^h\, ,  \qquad   \mathrm{at} \ \ y=h(x) \, . \label{usol3} 
\end{align}

\begin{figure}[h!]
    \begin{center}
        \includegraphics{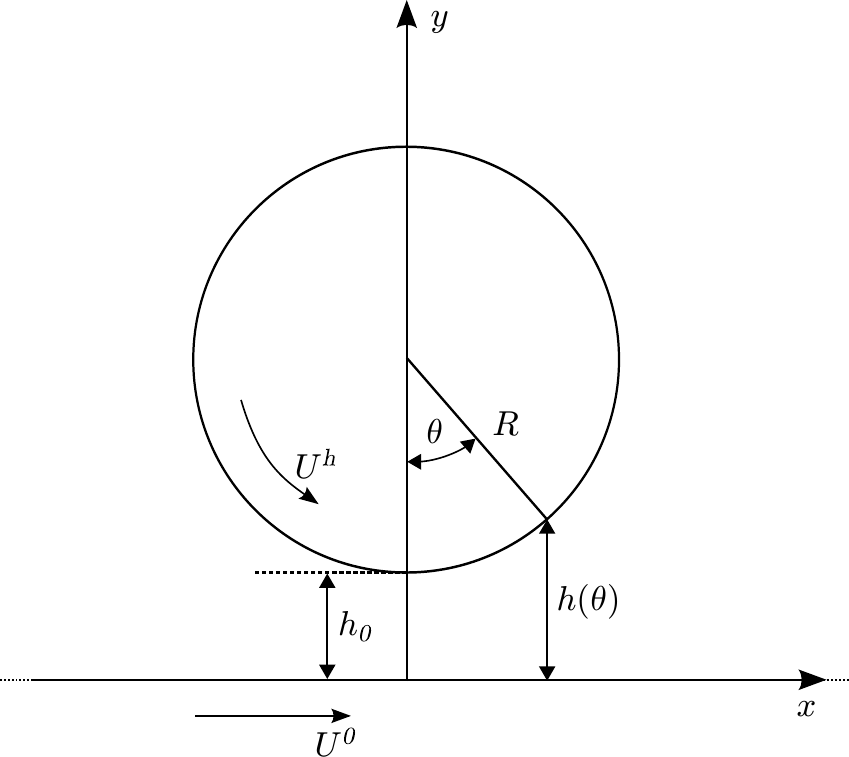}
    \end{center}
    \caption{The geometry of the cylinder-plane assembly.}
    \label{fig:cylinderplane}
\end{figure}

In order to illustrate the effect of the additional terms in the modified Reynolds equation,
we study the classic model problem of a rigid cylinder rolling on a plane.
(Refer to Szeri~\cite{Szeri2011} for an overview of the problem statement
and the traditional solution method.) Let $h_0$ be the minimum distance
between the cylinder of radius $R$ and the plane. Then the film thickness
$h(\theta)$ as a function of the angular coordinate $\theta$ is given by
\begin{equation}
    h(\theta) = -\frac{R}{n}(1+n \cos\theta),\quad n = - \frac{R}{h_0 + R}.
\end{equation}
The geometry of the problem is illustrated in Figure~\ref{fig:cylinderplane}.

After performing the change of variables $x = R \sin\theta$, the modified Reynolds equation \eqref{modified}
becomes
\begin{equation}
    \frac{\mathrm{d}}{\mathrm{d} \theta} \left[  \left( \frac{h^3}{12\mu }  - \mathcal{A}^2\mu  \, \int_0^{h} y(h-y)\, \left( \frac{\partial u^1}{\partial y} \right)^2 \, \mathrm{d} y \right) 
    \frac{\partial\theta}{\partial x}\frac{\mathrm{d} p^0}{\mathrm{d} \theta}    \right] = \frac{1}{2}(U^h-U^0)\frac{\partial h}{\partial \theta},
\label{modifiedtheta}
\end{equation}
where
\begin{equation}
    \frac{\partial \theta}{\partial x} = \frac{1}{R \cos \theta}, \quad
    \frac{\partial h}{\partial \theta} = R \sin \theta.
\end{equation}
Similarly, for equation \eqref{usol1} we get
\begin{equation}
    \mu \frac{\partial^2 u^1 }{\partial y^2}  +  \mathcal{A}^2\mu^2 \, \frac{\partial \theta}{\partial x} \frac{\partial p^0}{\partial \theta}  \left( \frac{\partial u^1}{\partial y} \right)^2 \, =  \frac{\partial \theta}{\partial x}\frac{\partial p^0}{\partial x}  \, . \label{usol1theta} 
\end{equation}

\renewcommand{\theenumi}{\alph{enumi}}

Equation \eqref{modifiedtheta} will be subject to the Swift--Stieber boundary conditions (see Szeri \cite{Szeri2011})
\begin{equation}
    \label{swiftstieber}
    p^0 = \begin{cases}
        0 & \text{when $\theta = -\frac{1}{2}\pi$},\\
        \frac{\mathrm{d}p^0}{\mathrm{d}\theta} = 0 & \text{when $\theta = \theta_2$},
    \end{cases}
\end{equation}
where yet another unknown $\theta_2$ is introduced through the unknown
location of the exit boundary. To solve the resulting system of nonlinear equations, we  perform iteratively the following
three-step process
\begin{enumerate}[Step 1.]
    \item \label{sa} Compute $(p^0,\theta_2)$ from \eqref{modifiedtheta}, \eqref{swiftstieber} with $u^1$ taken from the previous iteration (choose initially $u^1 = 0$).
    \item \label{sb} Substitute $p^0$  into \eqref{usol1theta} and solve for $u^1$ on a set of predefined angular
        values within the interval $[-\frac{1}{2}\pi,\theta_2]$.
    \item \label{sc} Interpolate $u^1$ between the chosen angular values and go to Step 1.
\end{enumerate}

\begin{figure}[h!]
    \begin{center}
        \includegraphics[width=0.6\textwidth]{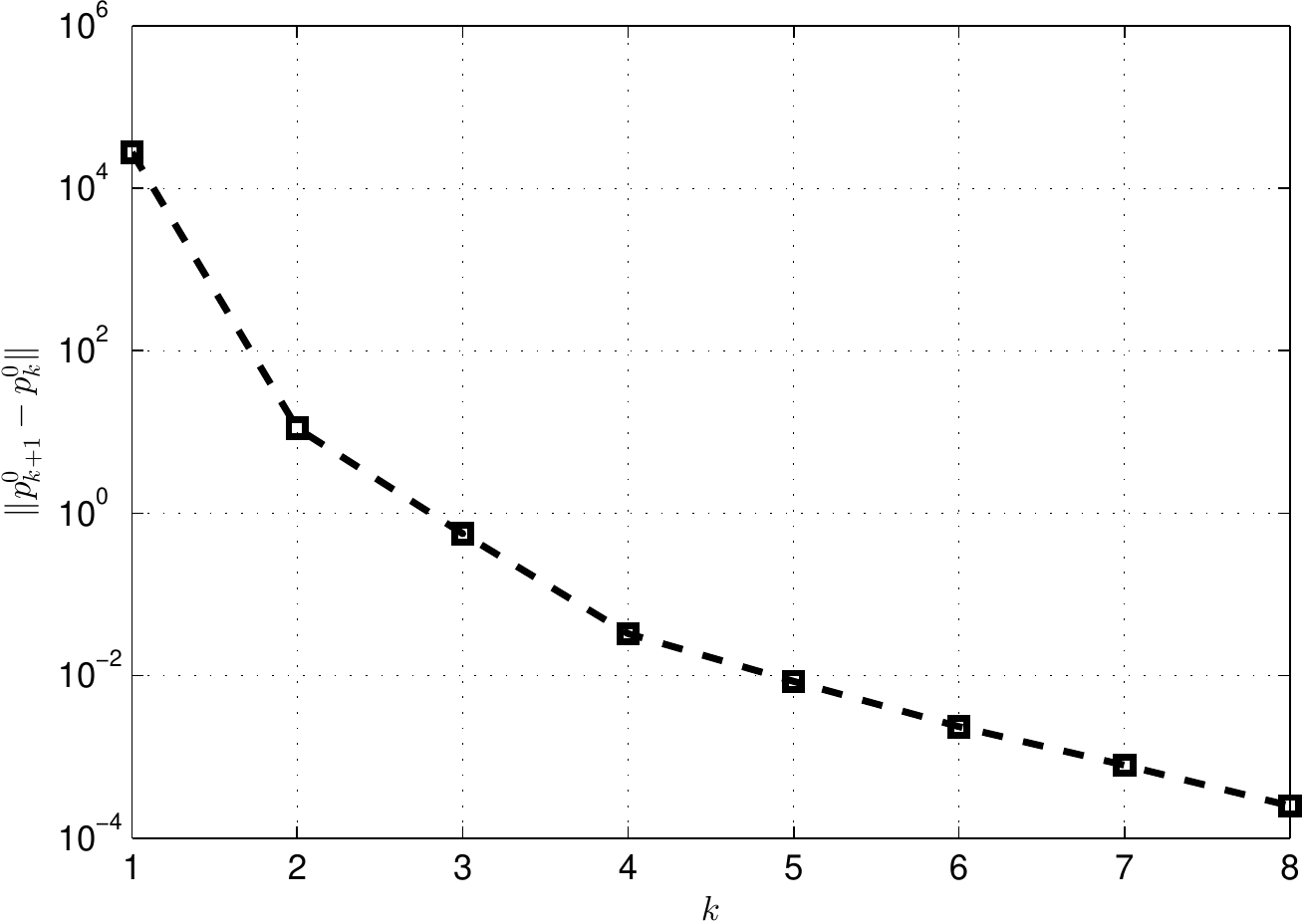}
    \end{center}
    \caption{Convergence of the difference of the pressure iterates in the $L^2$-norm with $\alpha = 5.59\cdot10^{-8}\,(\text{Pa})^{-1}$ and $h_0/R = 10^{-4}$. Here $p^0_k$ is the $k$'th iteration
       in the procedure  described in Steps \ref{sa}--\ref{sc} with $p^0_1$ corresponding to the
    solution of the modified Reynolds equation \eqref{modifiedtheta} when $u^1=0$.
   Given that the trailing edge angle $\theta_2$ might change from one iteration
    to another, the $L^2$-norm is computed using the current iterate's $\theta_2$ and
    the previous pressure is taken as zero  wherever the new $\theta_2$ is larger
    than the previous one.}
    \label{fig:convergence1}
\end{figure}

\begin{figure}[h!]
    \begin{center}
        \includegraphics[width=0.6\textwidth]{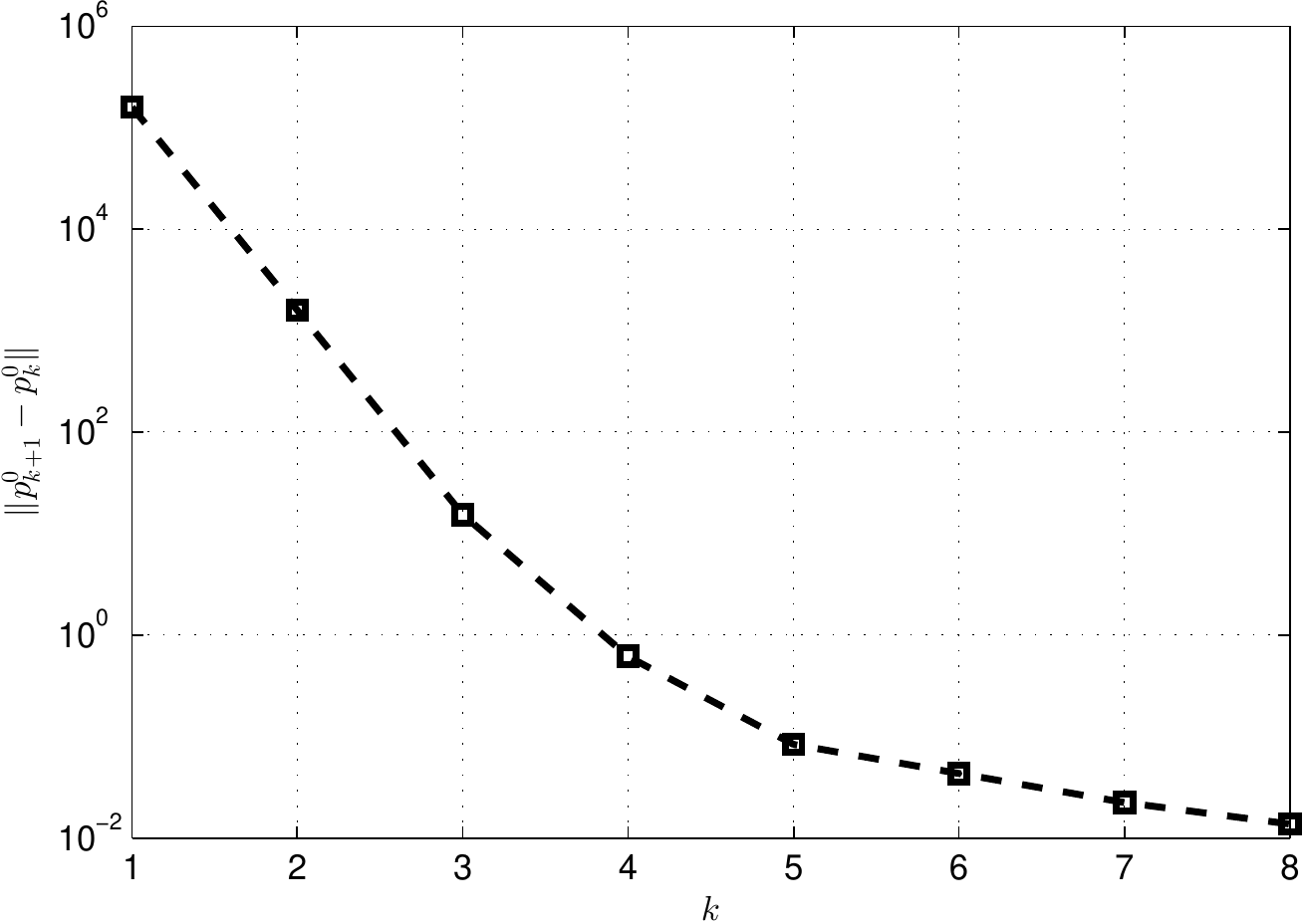}
    \end{center}
    \caption{Convergence of the  difference of the pressure iterates in the $L^2$-norm  when  $\alpha = 1.75\cdot10^{-7}\,(\text{Pa})^{-1}$ and $h_0/R = 10^{-3}$.}
    \label{fig:convergence2}
\end{figure}

The nonlinear equation  \eqref{modifiedtheta} at Step \ref{sa} is approximated
using a fixed-point iteration and the solution of the resulting linear equations
is done by the finite element method with
piecewise-linear elements. The finite element mesh is uniform in $\theta$
and $10^4$ elements were used. The correct value of $\theta_2$ at which
the condition \eqref{swiftstieber} is approximately satisfied is sought
through a binary search with initial lower and upper bounds at zero
and $\frac{1}{2}\pi$, respectively.

The values of $\theta$ at which  problem \eqref{usol1theta}, \eqref{usol2}, \eqref{usol3}
is solved in Step \ref{sb} are chosen to coincide with the nodes of
the finite element mesh of Step \ref{sa}. The nonlinear equation \eqref{usol1theta} 
is again approximated using a fixed-point scheme. The values of $\partial p^0/\partial \theta$ 
present in \eqref{usol1theta}
are computed from the solution obtained at the previous step. The
resulting linear equations are solved by the finite element
method with 50 uniformly distributed piecewise-linear elements.
Note that step (\ref{sb}) is perfectly parallelizable
for different values of $\theta$.

\begin{figure}[h!]
    \begin{center}
        \includegraphics[width=0.6\textwidth]{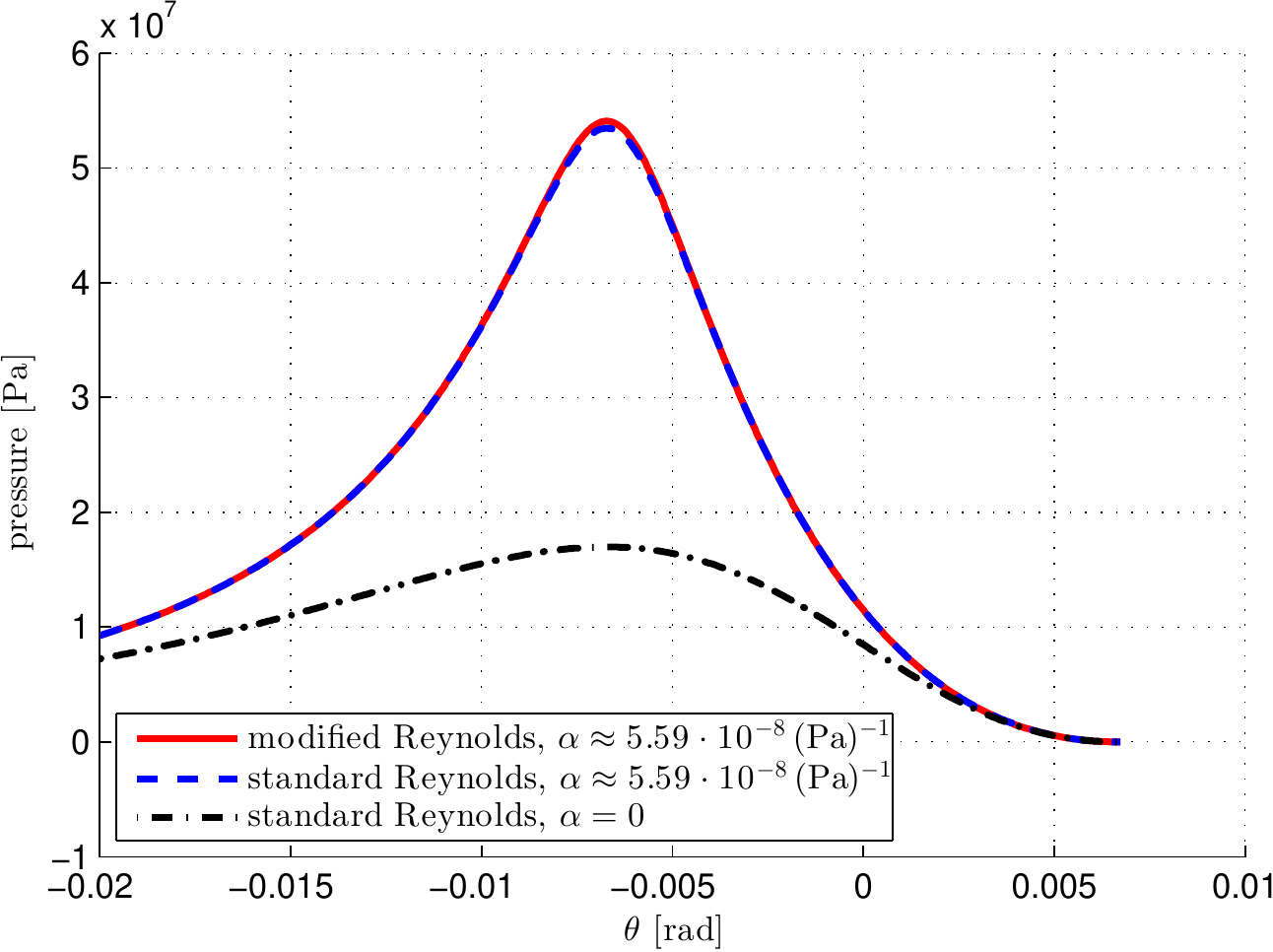}
    \end{center}
    \caption{The resulting pressure fields for the three models:
   standard Reynolds equation with a constant viscosity,  standard Reynolds equation with a pressure-dependent viscosity and  modified Reynolds equation with a pressure-dependent viscosity,
the last two plotted with $\alpha = 5.59 \cdot 10^{-8}\,(\text{Pa})^{-1}$
and $h_0/R=10^{-4}$.
Here $\theta_2\approx0.0067\,\text{rad}$ and the maximum pressure
values are $54.1\,\text{MPa}$ and $53.5\,\text{MPa}$ for the modified
and standard Reynolds equations.}
    \label{fig:pressure1}
\end{figure}

\begin{figure}[h!]
    \begin{center}
        \includegraphics[width=0.6\textwidth]{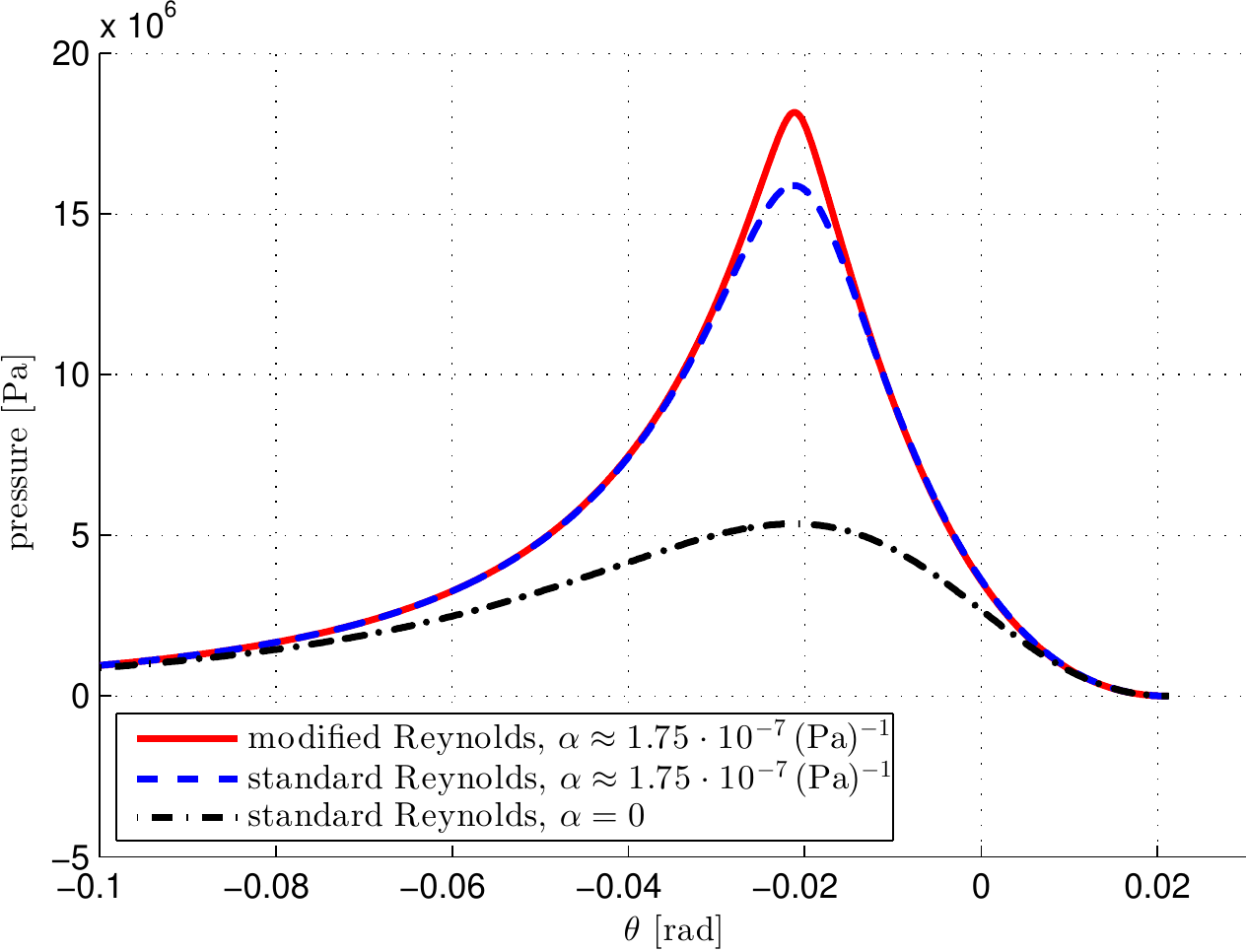}
    \end{center}
    \caption{The resulting pressure fields for the three models:
     standard Reynolds equation with a constant viscosity,  standard Reynolds equation with a pressure-dependent viscosity and  modified Reynolds equation with a pressure-dependent viscosity,
the last two plotted with $\alpha = 1.75 \cdot 10^{-7}\,(\text{Pa})^{-1}$
and $h_0/R=10^{-3}$.
Here $\theta_2\approx0.021\,\text{rad}$ and the maximum pressure
values are $18.2\,\text{MPa}$  for the modified
and $15.9\,\text{MPa}$ for the standard Reynolds equation.}
    \label{fig:pressure2}
\end{figure}

\begin{figure}[h!]
    \begin{center}
        \includegraphics[width=0.6\textwidth]{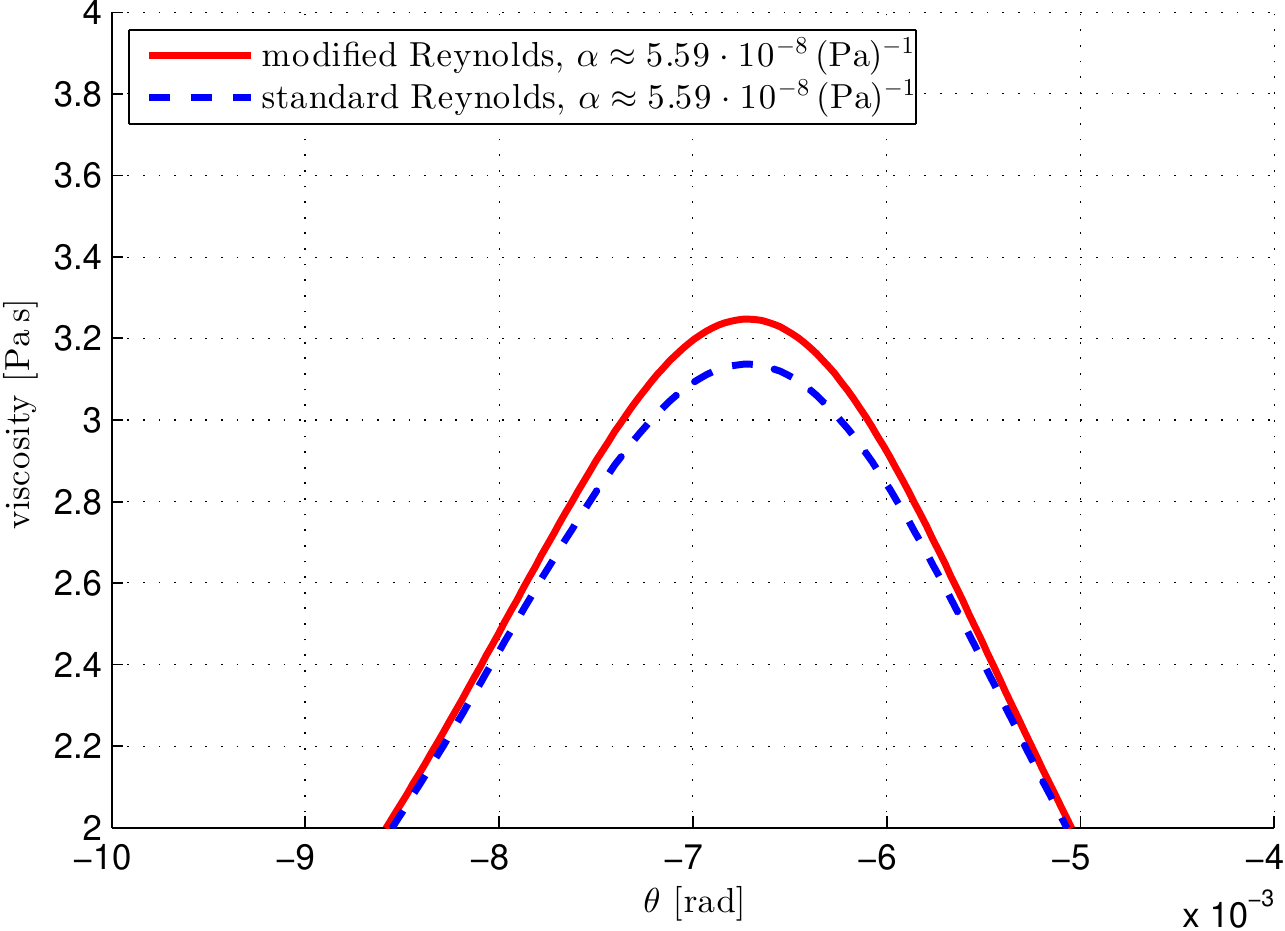}
    \end{center}
    \caption{The  viscosity distributions computed from the models with a pressure-dependent viscosity when $\alpha = 5.59 \cdot 10^{-8}\,(\text{Pa})^{-1}$
and $h_0/R=10^{-4}$. The maximum viscosity values
    are $3.24\,\text{Pa}\,\text{s}$  for the modified
and $3.14\,\text{Pa}\,\text{s}$ for the standard Reynolds equation. The corresponding pressure field is illustrated in Figure~\ref{fig:pressure1}.}
    \label{fig:viscosity1}
\end{figure}

A linear interpolation is performed in Step \ref{sc}
after which the newly solved $u^1$ is substituted
in \eqref{modifiedtheta}.
The convergence of the pressure between subsequent iterations
in $L^2$-norm with the parameter values $\alpha=5.59\cdot10^{-8}\,(\text{Pa})^{-1}$ and $\alpha=1.75\cdot10^{-7}\,(\text{Pa})^{-1}$,
$h_0/R = 10^{-4}$ and $h_0/R = 10^{-3}$,
$\mu_0 = 0.158\,\text{Pa}\,\text{s}$, $U^0=0\,\text{m}/\text{s}$ and $U^h=1\,\text{m}/\text{s}$ are visualized
in Figures~\ref{fig:convergence1} and \ref{fig:convergence2}.

\begin{figure}[h!]
    \begin{center}
        \includegraphics[width=0.6\textwidth]{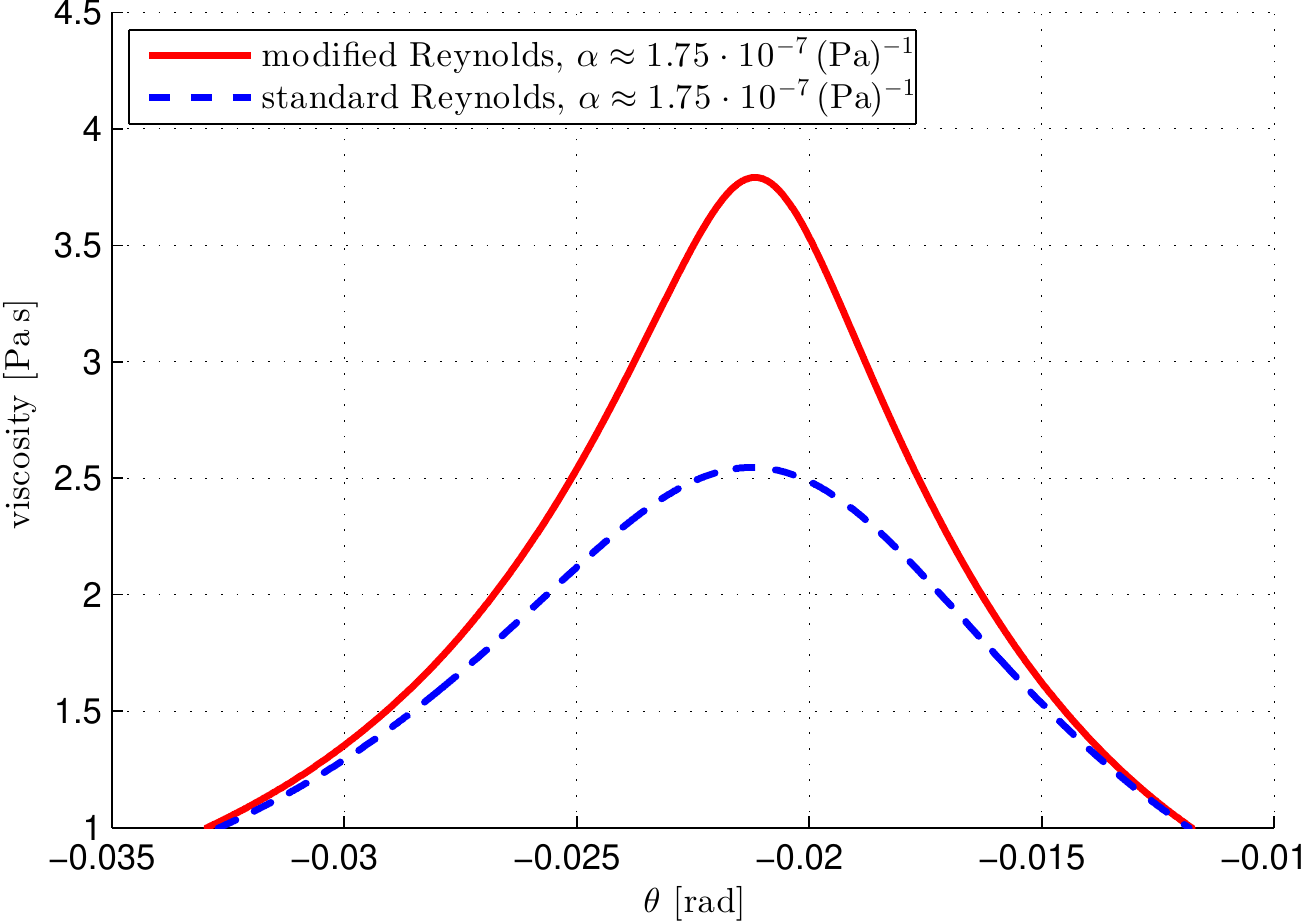}
    \end{center}
    \caption{The viscosity distributions computed from the  models with a pressure-dependent viscosity with $\alpha = 1.75 \cdot 10^{-7}\,(\text{Pa})^{-1}$
and $h_0/R=10^{-3}$. The maximum viscosity values
    are $3.8\,\text{Pa}\,\text{s}$  for the modified
and $2.5\,\text{Pa}\,\text{s}$ for the standard Reynolds equations. The corresponding pressure field is illustrated in Figure~\ref{fig:pressure2}.}
    \label{fig:viscosity2}
\end{figure}

The resulting pressure and viscosity fields after nine iterations are compared
to the solution of the standard Reynolds equation with a
pressure-dependent viscosity in Figures~\ref{fig:pressure1}, \ref{fig:pressure2}, \ref{fig:viscosity1} and \ref{fig:viscosity2}.
For comparison, the isoviscous case is also depicted   in Figures~\ref{fig:pressure1}, \ref{fig:pressure2}. 
The results reveal that
higher maximum pressure and viscosity values are obtained for
both sets of parameter values when
the modified Reynolds equation is applied.

\begin{figure}[h!]
    \begin{center}
        \includegraphics[width=0.7\textwidth]{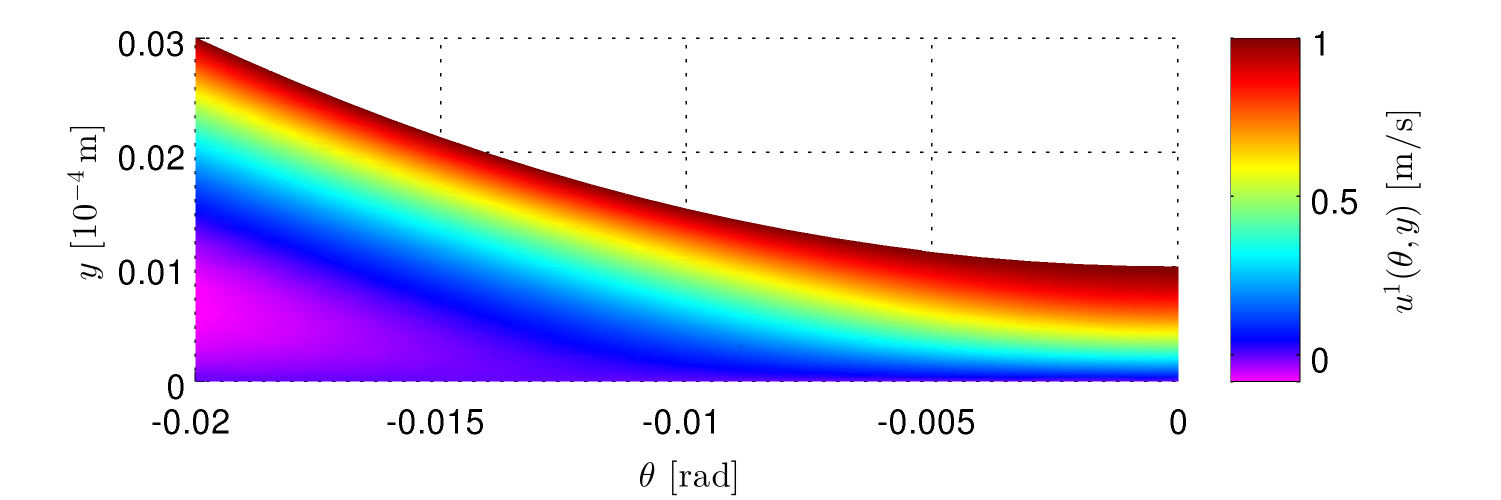}
        \includegraphics[width=0.7\textwidth]{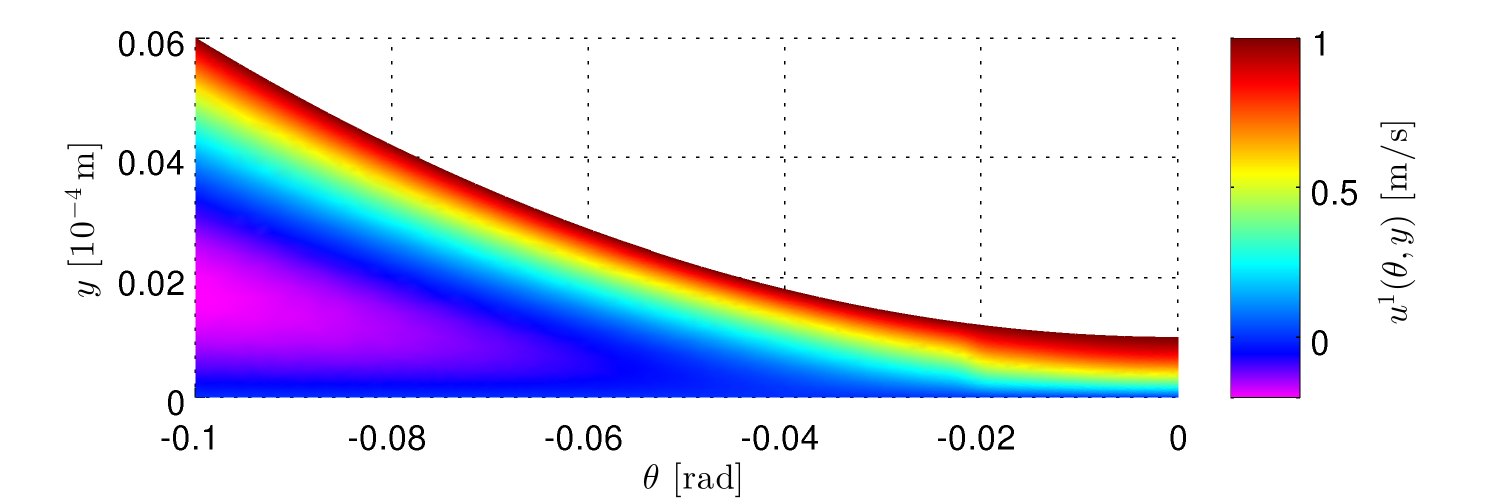}
    \end{center}
    \caption{The velocity field $u^1(\theta,y)$ corresponding to the pressure field of the final iteration with parameter values $\alpha = 5.59\cdot10^{-8}\,(\text{Pa})^{-1}$, $h_0/R = 10^{-4}$ (upper) and $\alpha = 1.75 \cdot 10^{-7}\,(\text{Pa})^{-1}$, $h_0/R = 10^{-3}$ (lower).
    Note that the lower left area in both drawings corresponds to negative velocity suggesting backflow of lubricant before the pressure spike.}
    \label{fig:velocity}
\end{figure}

The pressure and viscosity  differences predicted by the two models and shown in Figures~\ref{fig:pressure2} and \ref{fig:viscosity2} are larger than the ones  in Figures~\ref{fig:pressure1} and \ref{fig:viscosity1}.  
 In Figures~\ref{fig:pressure2} and \ref{fig:viscosity2}, a higher ratio $h_0/R = 10^{-3}$ allows us to use a higher
value of the pressure-viscosity coefficient which leads to lower pressures but, with our modified model, to a significantly higher viscosity. On the other hand,  larger values of $\alpha$, with a fixed $h_0/R$ ratio,
 caused numerical instabilities. The corresponding velocity fields $u^1$
 are visualized in Figure~\ref{fig:velocity}.

 The computations were performed using MATLAB R2014a on a HP ProLiant BL460c server blade
 with two 8-core Intel Xeon E5-2670 processors and 256 GB of RAM provided
 by Aalto University IT Services.

\section{Conclusions}

We have proposed a new Reynolds type lubrication approximation for fluids with pressure dependent viscosities. 
Starting with the full balance of linear momentum  equations, we have derived  a coupled set of dimensionally reduced equations governing the flow of piezoviscous fluids in hydrodynamical lubrication problems. 

 As shown,  both with a rigorous analysis and through numerical computations, the correction to the classical Reynolds equation can  be significant for certain values of the pressure-viscosity coefficient. The largest deviation from the classical lubrication approximation seems to occur for higher values of the pressure-viscosity coefficient and in less thin domains, at least if one adopts the  Barus formula for the pressure-viscosity relationship.  One should keep in mind though that the lower is the value of the pressure-viscosity coefficient and the thinner is the domain more concentrated, higher and less easily computable become the pressure spikes. 
 
Although derived for the Barus formula, our modified Reynolds approximation applies for more general pressure-viscosity relationship, e.g., the Roelands  formula. Its usefulness beyond the  study case presented here is under investigation.

\section*{Acknowledgements}
 
JHV was partially funded by the FCT Project PTDC/MAT-CAL/0749/2012. 
TG gratefully acknowledges the financial support from Tekes (Decision nr. 40205/12).


\section*{References}

\begin{thebibliography}{10}
\expandafter\ifx\csname url\endcsname\relax
  \def\url#1{\texttt{#1}}\fi
\expandafter\ifx\csname urlprefix\endcsname\relax\def\urlprefix{URL }\fi
\expandafter\ifx\csname href\endcsname\relax
  \def\href#1#2{#2} \def\path#1{#1}\fi

\bibitem{Rey1886}
O.~Reynolds, On the theory of lubrication and its application to {M}r {T}ower's
  experiments, Philosophical Transactions of the Royal Society of London 177
  (1886) 159--209.

\bibitem{Bridgman31}
P.~Bridgman, The Physics of High Pressure, MacMillan, 1931.

\bibitem{Bair2007}
S.~Bair, High Pressure Rheology for Quantitative Lubricants, Elsevier, 2007.

\bibitem{Szeri2011}
A.~Z. Szeri, Fluid Film Lubrication, 2nd ed., Cambridge University Press, 2011.

\bibitem{DH66}
D.~Dowson, G.~Higginson, Elasto-Hydrodynamic Lubrication: The Fundamentals of
  Roller and Gear Lubrication, Pergamon Press, 1966.

\bibitem{RS03}
{\relax K.R}.~Rajagopal, A.~Szeri, On an inconsistency in the derivation of the
  equations of elastohydrodynamic lubrication, Proceedings of the Royal Society
  of London. Series A. Mathematical Physical and Engineering Sciences 459
  (2003) 2771--2786.

\bibitem{Ren86}
M.~Renardy, Some remarks on the {N}avier--{S}tokes equations with a
  pressure-dependent viscosity, Communications in Partial Differential
  Equations 11 (1986) 779--793.

\bibitem{Gazz97}
F.~Gazzola, A note on the evolution {N}avier--{S}tokes equations with a
  pressure-dependent viscosity, Zeitschrift f\"{u}r angewandte Mathematik und
  Physik 48 (1997) 760--773.

\bibitem{HMR01}
J.~Hron, J.~M\'{a}lek, {\relax K.R}.~Rajagopal, Simple flows of fluids with
  pressure-dependent viscosities, Proceedings of the Royal Society of London.
  Series A. Mathematical Physical and Engineering Sciences 257 (2001)
  1603--1622.

\bibitem{MNR02}
J.~M\'{a}lek, J.~Ne\v{c}as, {\relax K.R}.~Rajagopal, Global analysis of the
  flows of fluids with pressure-dependent viscosities, Archive for Rational
  Mechanics and Analysis 165 (2002) 243--269.

\bibitem{FMR05}
M.~Franta, J.~M\'{a}lek, {\relax K.R}.~Rajagopal, On steady flows of fluids
  with pressure- and shear-dependent viscosities, Proceedings of the Royal
  Society of London. Series A. Mathematical Physical and Engineering Sciences
  461 (2005) 651--670.

\bibitem{MR07}
J.~M\'{a}lek, {\relax K.R}.~Rajagopal, Mathematical properties of the equations
  governing the flow of fluids with pressure and shear rate dependent
  viscosities, in: S. Friedlander and D. Serre (eds.), Handbook of Mathematical
  Fluid Dynamics, Vol.~4, Elsevier, 2007.

\bibitem{BMR09}
M.~Bul\'{i}\v{c}ek, J.~M\'{a}lek, {\relax K.R}.~Rajagopal, Mathematical
  analysis of unsteady flows of fluids with pressure, shear-rate and
  temperature dependent material moduli, that slip at solid boundaries, SIAM
  Journal of Mathematical Analysis 41 (2009) 665--707.

\bibitem{SV10}
G.~Saccomandi, L.~Vergori, Piezo-viscous flows over an inclined surface, The
  Quarterly of Applied Mathematics 68 (2010) 747--763.

\bibitem{HLS12}
A.~Hirn, M.~Lanzend\"{o}rder, J.~Stebel, Finite element approximation of flow
  of fluids with shear-rate- and pressure-dependent viscosity, IMA Journal of
  Numerical Analysis 32 (2012) 1604--1634.

\bibitem{BCGV13}
G.~Bayada, B.~Cid, G.~Garc\'{i}a, C.~V\'{a}zquez, A new more consistent
  {R}eynolds model for piezoviscous hydrodynamic lubrication problems in line
  contact devices, Applied Mathematical Modelling 37 (2013) 8505--8517.

\bibitem{NV07}
S.~Nazarow, J.~Videman, A modified nonlinear {R}eynolds equation for thin
  viscous flows in lubrication, Asymptotic Analysis 52 (2007) 1--36.

\bibitem{SGRW00}
C.~Schafer, P.~Giese, W.~Rowe, N.~Woolley, Elastohydrodynamically lubricated
  line contact based on the Navier-Stokes equations, in: Thinning Films and
  Tribological Interfaces, Proceedings of 26th Leeds--Lyon Symposium (ed. D.
  Dowson), Elsevier, 2000.

\bibitem{Barus93}
C.~Barus, Isothermals, isopiestics and isometrics relative to viscosity,
  American Journal of Science 45 (1893) 87--96.

\bibitem{Roe66}
C.~Roelands, Correlational aspects of the viscosity-temperature-pressure
  relationship of lubricating oils, Ph.D. thesis, Technische Hogeschool Delft,
  The Netherlands (1966).

\bibitem{PDR99}
M.~Paluch, Z.~Dendzik, S.~Rzoska, Scaling of high-pressure viscosity data in
  low-molecular-weight glass-forming liquids, Physical Review B 60 (1999)
  2979--2982.

\bibitem{BK03}
S.~Bair, P.~Kottke, Pressure-viscosity relationships for elastohydrodynamics,
  Tribology Transactions 46 (2003) 289--295.

\bibitem{BLW06}
S.~Bair, Y.~Liu, G.~J. Wang, The pressure-viscosity coefficient for {N}ewtonian
  {EHL} film thickness with general piezoviscous response, Journal of Tribology
  128 (2006) 624--631.

\bibitem{HSJ04}
B.~Hamrock, S.~Schmid, B.~Jacobson, Fundamentals of Fluid Film Lubrication, 2nd
  ed., Marcel Dekker, 2004.

\bibitem{vanLee09}
H.~van Leeuwen, The determination of the pressure-viscosity coefficient of a
  lubricant through an accurate film thickness formula and accurate film
  thickness measurements, Proceedings of the Institution of Mechanical
  Engineers, Part J: Journal of Engineering Tribology 212 (2009) 1143--1163.

\bibitem{PSR12}
V.~Pr\r{u}\v{s}a, S.~Srinivasan, {\relax K.R}.~Rajagopal, Role of pressure
  dependent viscosity in measurements with falling cylinder viscometer,
  International Journal of Non-linear Mechanics 47 (2012) 743--750.

\bibitem{Bur39}
J.~Burgers, Mechanical considerations, model systems, phenomenological theories
  of relaxation and viscosity, in: First Report on Viscosity and Plasticity,
  Nordemann Publishing, 1939.

\bibitem{Max66}
J.~Maxwell, On the dynamical theory of gases, Philosophical Transactions of the
  Royal Society A: Mathematical, Physical and Engineering Sciences 157 (1866)
  26--78.

\bibitem{Old50}
J.~Oldroyd, On the formulation of rheological equations of state, Proceedings
  of the Royal Society A: Mathematical, Physical and Engineering Sciences 200
  (1960) 523--591.

\end{thebibliography}
\end{document}